# A STUDY ON THE QUADRATIC JACO GRAPH

**R. Jayasree**

Research Scholar, Department of Mathematics

Nirmala College for Women

Coimbatore, India

sree.ramachandran14@gmail.com

**A. Kulandai Therese**

Department of Mathematics

Nirmala College for Women

Coimbatore, India

infanta1960@gmail.com

**U. Mary**

Department of Mathematics

Nirmala College for Women

Coimbatore, India

marycbe@gmail.com

**Johan Kok**

Tshwane Metropolitan Police Department

City of Tshwane, Republic of South Africa

Kokkiek2@tshwane.gov.za

**Abstract:**

In this paper we introduce the quadratic Jaco graph. The characteristics, properties and some graph invariants of quadratic Jaco graphs are also discussed. The observation that quadratic Jaco graphs are well-structured in respect of complete graphs and bridges led to the notion of braided complete graphs.

## 1. INTRODUCTION

For a general reference to notation and concepts of graph and digraphs we refer to [1].

In the early work the family of finite Jaco graphs (order 1) and (order *a*) were finite directed graphs derived from an infinite directed graph (order 1) and (order *a*) called the 1-root digraph and *a*-root digraph, respectively [4, 5]. Later the concept of linear Jaco graphs was introduced by Kok et al. [6] and it effectively redefined the initial concept. Linear Jaco graphs are a family of finite directed graphs which are derived from an infinite directed graph called the *f(x)*-root digraph. The incidence function is a linear function $f(x) = mx + c$, $x \in N$, $m, c \in N_0$. The *f(x)*-root digraph is denoted by, $J_\infty(f(x))$. We note that no classification for the special case $m = 0$ has been considered. For $m = 0$ and $c > 0$, the corresponding finite Jaco graph is the union (disjoint union) of complete digraphs. Also the components are of maximum order, $c + 1$. We will call these Jaco graphs, *constant* Jaco graphs denoted, *c*-Jaco graphs.

For $m > 0$, finite Jaco graphs are simple, connected and directed graphs. We will call these *connected linear* Jaco graphs denoted, *l*-Jaco graphs. Unless stated otherwise the *l*-Jaco graph will mean a finite, connected linear Jaco graph. This classification clarifies the distinction between other polynomial incidence functions such as quadratic functions which will be the main study area.

## 2. QUADRATIC JACO GRAPH

### 2.1 On *q*-Jaco graphs for $f(x) = ax^2 + bx + c$

In this section we report on properties, results and concepts related to the family of quadratic Jaco graphs, also called *q*-Jaco graphs. These digraphs have the incidence function $f(x) = ax^2 + bx + c$, with $a, x \in N$, $b, c \in N_0$.

**Definition 2.1.1**

The infinite *q*-Jaco graph or root *q*-Jaco graph $J_\infty(f(x))$, $x \in N$ is defined by $V(J_\infty(f(x))) = \{v_i : i \in N\}$, $A(J_\infty(f(x))) \subseteq \{(v_i, v_j) : i, j \in N, i < j\}$ and $(v_i, v_j) \in A(J_\infty(f(x)))$ if and only if $[ai^2 + (b+1)i + c] - d^-(v_i) \geq j$.

**Definition 2.1.2**

The family of finite $q$-Jaco graphs denoted by $\{J_n(f(x)): f(x) = ax^2 + bx + c;\ x, m \in N$ and $c \in N_0\}$ is defined by $V(J_n(f(x))) = \{v_i : i \in N, i \leq n\}$, $A(J_n(f(x))) \subseteq \{(v_i, v_j) : i, j \in N, i < j \leq n\}$ and $(v_i, v_j) \in A(J_n(f(x)))$ if and only if $[ai^2 + (b+1)i + c] - d^-(v_i) \geq j$.

**Definition 2.1.3**

Vertices with degree $\Delta(J_n(f(x)))$ are called Jaconian vertices and the set of vertices with maximum degree is called the Jaconian set of a $q$-Jaco graph $J_n(f(x))$ and denoted by $\mathbb{J}(J_n(f(x)))$.

**Definition 2.1.4**

The lowest numbered (indexed) Jaconian vertex is called the prime Jaconian vertex of a $q$-Jaco graph.

**Definition 2.1.5**

If $v_i$ is the prime Jaconian vertex the complete subgraph on vertices $v_{i+1}, v_{i+2}, ..., v_n$ is called the Hope subgraph of a $q$-Jaco graph.

**2.2 Basic Properties of $q$-Jacograph**

**Property 2.2.1**

From the definition of a $q$-Jaco graph $J_n(f(x))$ it follows that, if for the prime Jaconian vertex $v_i$, we have $d(v_i) = f(i)$, then in the underlying $q$-Jaco graph denoted $J_n^*(f(x))$ we have $d(v_m) = f(m)\ \forall m \in \{1,2,3,...,i\}$.

**Property 2.2.2**

From the definition of a $q$-Jaco graph $J_n(f(x))$, it follows that $\Delta(J_k(f(x))) \leq \Delta(J_n(f(x)))$ $\forall k \leq n$.

**Property 2.2.3**

From the definition of a $q$-Jaco graph $J_n(f(x))$, it follows that the lowest degree attained by all Jaco graphs is $0 \leq \delta(J_n(f(x))) \leq f(1)$.

**Property 2.2.4**

The $d^-(v_k)$ for any vertex $v_k$ of a $q$-Jaco graph $J_n(f(x))$, $n \geq k$ is equal to $d(v_k)$ in the underlying Jaco graph $J_k^*(f(x))$.

### 2.3 Results on *q*-Jaco graphs

Our first result addresses some distinction between *constant, linear and quadratic* Jaco graphs.

**Lemma 2.3.1**

A complete graph $K_n, n \geq 1$ can be described by the underlying graph of some finite Jaco graph, $J_n(f(x))$ such as:

(i) $f(x) = c, c \geq n-1$ (*c*-Jaco graph) or,

(ii) $f(x) = mx + c, m \in N, c \in N_0$ with $m + c \geq n - 1$ (*l*-Jaco graph) or,

(iii) $f(x) = ax^2 + bx + c, a \in N, b, c \in N_0$ with $a + b + c \geq n - 1$ (*q*-Jaco graph).

**Proof:**

Let us consider $f(x) = c, c \geq n - 1$.

(i) Clearly for $c = 0$ the complete graph $K_1$ is described. For $f(1) = c \geq 1$, the arc $(v_1, v_{c+1})$ exists, as the maximum reach of vertex $v_1$ in accordance with Definition 2.1.1. Therefore all arcs $(v_1, v_j), 2 \leq j \leq c+1$ exist. Similarly all arcs $(v_i, v_j), 2 \leq i \leq c$ and $i+1 \leq j \leq c+1$ exist. Clearly, the underlying finite constant Jaco graph is the complete graph $K_{c+1}$. Hence, for $J_l(f(x)), 1 \leq l \leq c+1$ the underlying finite constant Jaco graphs are complete graphs.

Now consider $f(x) = mx + c, m \in N, c \in N_0$ with $m + c \geq n - 1$.

(ii) Since $f(1) = m + c$, the arc $(v_1, v_{m+c+1})$ exists as the maximum reach of vertex $v_1$ in accordance with Definition 2.2.1. Therefore all arcs $(v_1, v_j)$, $2 \leq j \leq m+c+$ exist. Similarly all arcs $(v_i, v_j)$, $2 \leq i \leq m+c$ and $i+1 \leq j \leq m+c+1$ exist. Clearly, the underlying finite linear Jaco graph is the complete graph $K_{m+c+1}$. Hence, for $J_l(f(x)), 1 \leq l \leq m+c+1$, the underlying finite linear Jaco graphs are complete graphs.

Finally consider $f(x) = ax^2 + bx + c, m \in N, c \in N_0$ with $a + b + c \geq n - 1$.

(iii) Since $f(1) = a+b+c$, the arc $(v_1, v_{a+b+c+1})$ exists as the maximum reach of vertex $v_1$ in accordance with Definition 2.1.1. Therefore all arcs $(v_1, v_j)$, $2 \leq j \leq a+b+c+1$ exist.

Similarly all arcs $(v_i, v_j)$, $2 \leq i \leq a+b+c$ and $i+1 \leq j \leq a+b+c+1$ exist. Clearly, the finite underlying quadratic Jaco graph is the complete graph $K_{a+b+c+1}$. Hence, for $J_l(f(x))$, $1 \leq l \leq a+b+c+1$ the underlying finite quadratic Jaco graphs are complete graphs.

Hence, a complete graph $K_n$, $n \geq 1$ can be described by some finite underlying Jaco Graph, $J_n(f(x))$ where $f(x) = c, c \geq n-1$ or $f(x) = mx+c, m \in N, c \in N_0$

or $f(x) = ax^2 + bx + c, a \in N, b, c \in N_0$. $\square$

**Lemma 2.3.2**

For the q-Jaco graphs $J_i(f(x)), i \in \{1,2,3,..., f(1)+1\}$ we have $\Delta(J_i(f(x))) = i-1$ and $\mathcal{J}(J_i(f(x))) = \{v_k : 1 \leq k \leq i\} = V(J_i(f(x)))$.

**Proof:**

Clearly for $i = 1$ the corresponding q-Jaco graph is an isolated vertex $v_1$ or put differently, the complete graph $K_1$. Therefore the result holds for $i = 1$.

For any $2 \leq i \leq f(1)+1$ it follows from Definition 3.1.2 that the arcs $(v_1, v_2), (v_1, v_3), ..., (v_1, v_i)$ and $(v_2, v_3), (v_2, v_4), ..., (v_2, v_i)$ and $(v_3, v_4), (v_3, v_5), ..., (v_3, v_i)$ and ... and $(v_{i-1}, v_i)$ exist. Hence, the underlying graph is a complete graph and all vertices $v_k$, $1 \leq k \leq i$ have degree $i - 1$. So in terms of the corresponding q-Jaco graph we have $\mathcal{J}(J_i(f(x))) = \{v_k : 1 \leq k \leq i\} = V(J_i(f(x)))$. Finally since no arc $(v_1, v_i), i > f(1) + 1$ is defined the range $i \in \{1,2,3,..., f(1)+1\}$ holds. $\square$

**Proposition 2.3.3**

For the q-Jaco graphs $J_{f(i)+1}(f(x)), 1 \leq i \leq f(1)+1$ the prime Jaconian vertex is the vertex $v_i$ and $\mathcal{J}(J_{f(i)}(f(x))) = \{v_l : i \leq l \leq f(1)+1\}$.

**Proof:**

From Definition 2.1.2 it follows that $J_{f(1)+1}(f(x))$ is a complete digraph $K_{f(1)+1}$. Therefore,

$\mathcal{G}(J_{f(1)+1}(f(x))) = \{v_l : 1 \le l \le f(1)+1\}$. Now for any $2 \le i \le f(1)+1$ it follows from Definition 2.1.2 that $d(v_i) = d(v_{i+1}) = \cdots = d(v_{f(1)+1}) = \Delta\left(J_{f(i)+1}(f(x))\right)$ and $d(v_{f(1)+1}) > d(v_j) \ge d(v_{j+1}) \ge d(v_{j+1}) \ldots \ge d(v_{f(i)+1}), f(i) + 2 \le j \le f(i) + 1$.

Hence, $v_i$ is the prime Jaconian vertex of $J_{f(1)+1}(f(x))$ and

$\mathcal{G}(J_{f(1)+1}(f(x))) = \{v_l : 1 \le l \le f(1)+1\}$. □

**Lemma 2.3.4**

If in a q-Jaco graph $J_n(f(x))$, and for the smallest $i$, we have $d(v_i) = f(i)$ and the arc $(v_i, v_n)$ is defined, then $v_i$ is the prime Jaconian vertex of $J_n(f(x))$.

**Proof:**

If in the construction of the q-Jaco graph $J_n(f(x))$ and for the smallest $i$ we have $d(v_i) = f(i)$ and the arc $(v_i, v_n)$ is defined, then $d(v_j) = f(j) < d(v_i), 1 \le j < i$. Also, $d(v_k) \le d(v_i), i < k \le n$. Therefore, $v_i$ is the prime Jaconian vertex of $J_n(f(x))$. □

Observe that $\Delta(J_i(f(x)))$ and $d^-(v_i)$ might repeat as $i$ increases to $i + 1$ and on an increase, the increase is always +1.

**Lemma 2.3.6**

For all q-Jaco graphs, $J_n(f(x)), n \ge 2$ and, $v_i, v_{i-1} \in V(J_n(f(x)))$ we have in the *underlying q-Jaco graph* $J_n^*(f(x))$, that $|d(v_i) - d(v_{i-1})| \le a(2i-1) + b$.

**Proof:**

Clearly, $\max|d(v_i) - d(v_{i-1})| = \max|f(i) - f(i-1)| = f(i) - f(i-1)$

$= ai^2 + bi + c - (a(i-1)^2 + b(i-1) + c)$

$= ai^2 + bi + c - ai^2 + 2ai - a - bi + b - c$

$= a(2i - 1) + b.$

Hence, $|d(v_i) - d(v_{i-1})| \le a(2i-1) + b.$ □

**Theorem 2.3.7**

The q-Jaco Graph $J_k(f(x)), k = f(f(1)) - f(1) + 1,$ is the smallest q-Jaco graph in $\{J_n(f(x)) : n \in N\}$ which has $\Delta(J_k(f(x))) = f(f(1))$ with the prime Jaconian vertex $v_{(a+b+c)}$.

**Proof:**

Clearly $d^-(v_{f(1)}) = f(1) - 1$. Furthermore, the maximum defined degree for vertex $v_{f(1)}$ is $f(f(1))$. Hence, at vertex $v_{f(f(1))-f(1)+1}$ in a sufficient large finite $q$-Jaco graph the defined out-arcs of vertex $v_{f(1)}$ including the arc $(v_{f(1)}, v_{f(f(1))-f(1)+1})$ ensure $d(v_{f(1)}) = f(f(1))$. Hence, the $q$-Jaco graph $J_k(f(x)), k = f(f(1)) - f(1) + 1$ is the smallest $q$-Jaco graph having $d(v_{f(1)}) = f(f(1))$.

Now since $d(v_i) < f(f(1)), 1 \leq i \leq f(1) - 1$ and $d(v_j) \leq f(f(1)), f(1) + 1 \leq j \leq f(f(1)) - f(1) + 1$ we have $\Delta(J_k(f(x))) = f(f(1))$ with the prime Jaconian vertex $v_{(a+b+c)}$ for $J_k(f(x)), k = f(f(1)) - f(1) + 1$. $\square$

## 2.4 $q$-Jaco graphs for $f(x) = x^2$

Unless mention otherwise, we shall in the main limit this section to the case $f(x) = x^2$.

Table 1 shows the results for the application of the Fisher algorithm, for $i \leq 35$.

Table 1

| $\phi(v_i) \to i \in N$ | $d^-(v_i)$ | $d^+(v_i) = i - d^-(v_i)$ | $\mathfrak{J}(J_i(x^2))$ | $\Delta(J_i(x^2))$ | $d_{J_i(x^2)}(v_1, v_i)$ |
|---|---|---|---|---|---|
| 1 | 0 | 1 | $\{v_1\}$ | 0 | 0 |
| 2 | 1 | 3 | $\{v_1, v_2\}$ | 1 | 1 |
| 3 | 1 | 8 | $\{v_2\}$ | 2 | 2 |
| 4 | 2 | 14 | $\{v_2\}$ | 3 | 2 |
| 5 | 3 | 22 | $\{v_2\}$ | 4 | 2 |
| 6 | 3 | 33 | $\{v_2, v_3, v_4, v_5\}$ | 4 | 3 |
| 7 | 4 | 45 | $\{v_3, v_4, v_5\}$ | 5 | 3 |
| 8 | 5 | 59 | $\{v_3, v_4, v_5\}$ | 6 | 3 |
| 9 | 6 | 73 | $\{v_3, v_4, v_5\}$ | 7 | 3 |
| 10 | 7 | 93 | $\{v_3, v_4, v_5\}$ | 8 | 3 |
| 11 | 8 | 113 | $\{v_3, v_4, v_5\}$ | 9 | 3 |
| 12 | 8 | 136 | $\{v_4, v_5\}$ | 10 | 4 |
| 13 | 9 | 160 | $\{v_4, v_5\}$ | 11 | 4 |
| 14 | 10 | 186 | $\{v_4, v_5\}$ | 12 | 4 |
| 15 | 11 | 214 | $\{v_4, v_5\}$ | 13 | 4 |
| 16 | 12 | 244 | $\{v_4, v_5\}$ | 14 | 4 |
| 17 | 13 | 276 | $\{v_4, v_5\}$ | 15 | 4 |
| 18 | 14 | 310 | $\{v_4, v_5\}$ | 16 | 4 |
| 19 | 14 | 347 | $\{v_5\}$ | 17 | 5 |
| 20 | 15 | 385 | $\{v_5\}$ | 18 | 5 |
| 21 | 16 | 425 | $\{v_5\}$ | 19 | 5 |
| 22 | 17 | 467 | $\{v_5\}$ | 20 | 5 |
| 23 | 18 | 511 | $\{v_5\}$ | 21 | 5 |
| 24 | 19 | 557 | $\{v_5\}$ | 22 | 5 |
| 25 | 20 | 605 | $\{v_5\}$ | 23 | 5 |
| 26 | 21 | 655 | $\{v_5\}$ | 24 | 5 |
| 27 | 22 | 707 | $\{v_5\}$ | 25 | 5 |
| 28 | 22 | 762 | $\{v_5, v_6, v_7, v_8, v_9, v_{10}, v_{11}\}$ | 25 | 6 |
| 29 | 23 | 818 | $\{v_6, v_7, v_8, v_9, v_{10}, v_{11}\}$ | 26 | 6 |
| 30 | 24 | 876 | $\{v_6, v_7, v_8, v_9, v_{10}, v_{11}\}$ | 27 | 6 |
| 31 | 25 | 959 | $\{v_6, v_7, v_8, v_9, v_{10}, v_{11}\}$ | 28 | 6 |
| 32 | 26 | 998 | $\{v_6, v_7, v_8, v_9, v_{10}, v_{11}\}$ | 29 | 6 |
| 33 | 27 | 1062 | $\{v_6, v_7, v_8, v_9, v_{10}, v_{11}\}$ | 30 | 6 |
| 34 | 28 | 1128 | $\{v_6, v_7, v_8, v_9, v_{10}, v_{11}\}$ | 31 | 6 |
| 35 | 29 | 1196 | $\{v_6, v_7, v_8, v_9, v_{10}, v_{11}\}$ | 32 | 6 |

**Proposition 2.4.1**

If for $n \in N$ the in-degree $d^-(v_n) = l = d^-(v_{n+1})$ then, $|\mathcal{J}(J_n(x^2))| \neq |\mathcal{J}(J_{n+1}(x^2))|$.

**Proof:**

Assume that for $n \in N$ the in-degree $d^-(v_n) = l = d^-(v_{n+1})$. We consider two cases.

Case 1. Let $|\mathcal{J}(J_n(x^2))| \geq 2$. Let vertex $v_i$ be the prime Jaconian vertex of $J_n(x^2)$. Now expand to the Jaco graph $J_{n+1}(x^2)$. Since, $d^-(v_n) = l = d^-(v_{n+1})$ the arc$(v_i, v_{n+1})$ cannot exist. Therefore, $|\mathcal{J}(J_n(x^2))| \neq |\mathcal{J}(J_{n+1}(x^2))|$.

Case 2. Let $|\mathcal{J}(J_n(x^2))| = 1$. Let vertex $v_i$ be the unique prime Jaconian vertex of $J_n(x^2)$. It implies that $d(v_i) > d(v_{i+1}) \geq d(v_{i+2}) \geq \cdots \geq d(v_k)$ for some $k$. It also implies that $d(v_i) = d(v_{i+1}) + 1$.

Expand to the Jaco graph $J_{n+1}(x^2)$. Since, $d^-(v_n) = l = d^-(v_{n+1})$ the arc$(v_i, v_{n+1})$ cannot exist. This implies that in $J_{n+1}(x^2)$, $d(v_i) = d(v_{i+1}) \geq d(v_{i+2}) \geq \cdots \geq d(v_k)$ for some $k$. Therefore, $|\mathcal{J}(J_{n+1}(x^2))| \geq 2$. Hence, $|\mathcal{J}(J_n(x^2))| \neq |\mathcal{J}(J_{n+1}(x^2))|$. □

**Lemma 2.4.2**

The in-degrees of the $q$-Jaco graphs for the incidence function $f(x) = x^2$ are stepwise consecutive non-negative integers.

**Proof:**

The result holds as the $q$-Jaco graph expands from order $i = 1$ to $i = 2$. Assume it holds for the $q$-Jaco graphs as the order expands from $i = k - 1$ to $i = k$. It means that $d^-(v_k) \in \{d^-(v_{k-1}), d^-(v_{k-1}) + 1\}$.

Now let the $q$-Jaco graph expand to order $i = k + 1$. Clearly, $N^-(v_{k+1}) = N^-(v_k) \cup \{v_k\}$ or $(N^-(v_k) - v_i) \cup \{v_k\}$, with $i$ the minimum for which arc $(v_i, v_k)$ exists. Whichever the case, it implies that $d^-(v_{k+1}) \in \{d^-(v_k), d^-(v_k) + 1\}$. Hence, the result holds for $q$-Jaco graphs of any order $i \in N$. □

We now state a general result which follows from Lemma 2.4.2.

**Corollary 2.4.3**

The result in Lemma 2.4.2 holds for all $q$-Jaco graphs with the incidence of the form $f(x) = ax^2 + bx + c$, with $a, x \in N$, $b, c \in N_0$.

**Proof:**

The proof is similar to that of Lemma 2.4.2. □

For the root $l$-Jaco graph we observed occasional repetition of the out-degree $d^+(v_i)$. We show that such repetition does not occur in the root $q$-Jaco graph for $f(x) = x^2$.

**Proposition 2.4.4**

Consider a $q$-Jaco graph $J_\infty(x^2)$, $n, x \in N$ then $d^+(v_i) \neq d^+(v_{i+1}), \forall i \in N$.

**Proof:**

From the previous table it follows that the result holds for $i = 1,2$. Assume the result holds for $1 \leq l \leq k$. Now consider $i = k + 1$.

Case 1. Let $d^-(v_k) = d^-(v_{k+1})$. Since, $k^2 - d^-(v_k) = k^2 - d^-(v_{k+1}) < k^2 - d^-(v_{k+1}) + (2k + 1) = (k + 1)^2 - d^-(v_{k+1})$, it follows that $d^+(v_k) \neq d^+(v_{k+1})$. Hence, the result holds for this case.

Case 2. Let $d^-(v_k) \neq d^-(v_{k+1})$. Then from Lemma 2.4.2 we have that $d^-(v_{k+1}) = d^-(v_k) + 1$. Therefore $k^2 - d^-(v_k) < k^2 - d^-(v_k) - 1 + 2k + 1 = (k + 1)^2 - (d^-(v_k) + 1) = (k + 1)^2 - d^-(v_{k+1})$. Hence, the result holds for this case.

$\therefore d^+(v_i) \neq d^+(v_{i+1}), \forall i \in N$. □

We now state a general result which follows from Proposition 2.4.4.

**Corollary 2.4.5**

The result in Proposition 2.4.4 holds for all $q$-Jaco graphs with the incidence of the form $f(x) = ax^2 + bx + c$, with $a, x \in N$, $b, c \in N_0$.

**Proof:**

The proof is similar to that of Proposition 2.4.4. □

## 2.5 Minimum and Maximum Chromatic sum of $q$-Jaco graphs

We begin with the concept of chromatic sums which were introduced by Kok et al. [8]. For concepts and notation not defined or explained in this section see [2, 7, 8].

### 2.5.1 The $\chi'$-Chromatic Sum and $\chi^+$-Chromatic Sum of Jaco Graphs [8]

If the colours represent different technology types and the configuration requirement is that at least one unit per technology type must be placed at a point in a network without similar technology types being adjacent, two further considerations come into play. Firstly, if the higher indexed colours represent technology types with higher failure rate (risk) then the placement of the maximal number of higher indexed colours units is the solution to ensure a functional network. On the other hand, if the lower indexed colours represent a less costly (procurement, installation, commissioning and maintenance) technology type, and minimising total cost is the priority, then the placement of maximal number of lower indexed units is the desired solution.

Let $C = \{c_1, c_2, c_3, ..., c_k\}$ allow a colouring $S$ of G. As stated in [2] there are $k!$ ways of allocating the colours to the vertices of G. Let the colour weight $\theta(c_i)$ be the number of times a colour $c_i$ is allocated to vertices. In general we refer to the colour sum of a colouring $S$ and define it,

$$\omega(S) = \sum_{i=1}^{k} i.\theta(c_i).$$

An interesting new invariant in respect of graph colouring was introduced in [8].

**Definition 2.5.2**

For a graph G the $\chi'$-Chromatic sum is defined to be:

$$\chi'(G) = \min\{\sum_{i=1}^{k} i\,\theta(c_i) : \forall\ \min\ proper\ colourings\ of\ G\}.$$

Such a minimum colouring sum is obtain by a greedy algorithm, i.e. colour the maximum number of vertices, allowed by the definition of a proper colouring, the colour $c_1$. This is followed by colouring the maximum number of the uncoloured vertices allowed by the definition of a proper colouring, the colour $c_2$, and so on. Iteratively complete the colouring of the graph. The colouring obtained is a proper colouring corresponding to the minimum colouring sum.

**Example 1:**

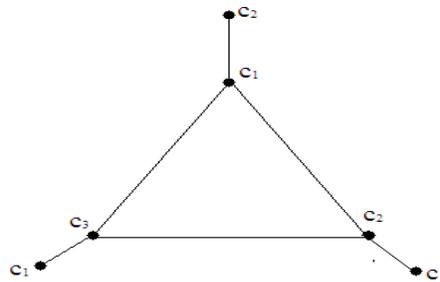

Fig 3.1: Minimum Chromatic Sum

**Definition 2.5.3**

For a graph G the $\chi^+$-Chromatic sum is defined to be:

$$\chi'(G) = \max\{\sum_{i=1}^{k} i\,\theta(c_i) : \forall \min \text{ proper colourings of } G\}.$$

A maximum colouring sum is obtain by considering the colour set $C=\{c_1,c_2,c_3,...,c_k\}$ which allows a minimum colouring sum and then, to recolour the vertices according to the colour mapping $c_i \mapsto c_{k-i+1}$.

**Example 2:**

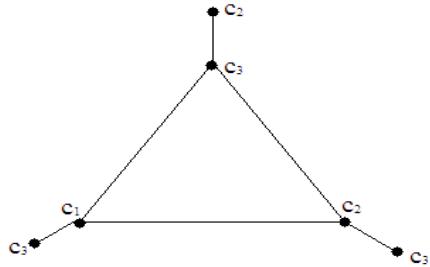

Fig 3.2 : Maximum Chromatic Sum

**Theorem 2.5.4**

For a complete graph $K_n$, the $\chi'$-chromatic sum and $\chi^+$-chromatic sum are given by

$$\chi'(K_n) = \chi^+(K_n) = \frac{n(n+1)}{2}.$$

**Proof**:

Let $C$ be a proper colouring of the complete graph. Then $C$ must contain at least $n$ colors say, $c_1, c_2, c_3, ..., c_n$. Because every vertex is adjacent to every other vertex in complete graph we have:

$\theta(c_1) = \theta(c_2) = ... = \theta(c_n) = 1,$ meaning that, $\theta(c_i) = 1, \forall c_i, 1 \leq i \leq n.$

Therefore, $1(1) + 2(1) + 3(1) + ... + n(1) \Rightarrow \chi'(K_n) = \frac{n(n+1)}{2}.$

Similarly, $\chi^+(K_n) = \frac{n(n+1)}{2}$. Hence, $\chi'(K_n) = \chi^+(K_n) = \frac{n(n+1)}{2}.$ □

**2.6  Arithmetic Mean and Variance of Chromatic Parameters**

Colouring the vertices of a graph G can be considered as a random experiment and a discrete variable X can be defined. Correspondingly the two major statistical (probability) parameters, the arithmetic mean and variance can be determined for the minimum and the maximum chromatic sums. This new concept was introduced by Sudev et al. in [9]. Following the probability mass function (pmf) the two probability parameters are defined as follows.

$$f(i) = \begin{cases} \dfrac{\theta(c_i)}{|V(G)|}, & i = 1,2,3,...,k, \\ 0 & elsewhere, \end{cases}$$

$$\text{and } \mu_C(G) = \dfrac{\sum_{i=1}^{k} i\,\theta(c_i)}{\sum_{i=1}^{k} \theta(c_i)} \text{ and } \sigma_C^2(G) = \dfrac{\sum_{i=1}^{k} i^2 \theta(c_i)}{\sum_{i=1}^{k} \theta(c_i)} - \left(\dfrac{\sum_{i=1}^{k} i\,\theta(c_i)}{\sum_{i=1}^{k} \theta(c_i)}\right)^2.$$

**Lemma 2.6.1**

The $\chi'$-chromatic mean, the $\chi^+$-chromatic mean, the $\chi'$-chromatic variance and the $\chi^+$-chromatic variance of a complete graph $K_n, n \geq 1$, are:

$$\mu_{\chi^+}(K_n) = \mu_{\chi'}(K_n) = \dfrac{n+1}{2} \text{ and } \sigma_{\chi^+}^2(K_n) = \sigma_{\chi'}^2(K_n) = \dfrac{n^2-1}{12}.$$

**Proof:**

For a complete graph $K_n$ a proper colouring is allowed by $n$ colours $\{c_1, c_2, c_3, ..., c_n\}$. Therefore, $\theta(c_i) = 1, \forall c_i, 1 \leq i \leq n$. By symmetry argument we have: $\chi(K_n) = \chi'(K_n) = \chi^+(K_n)$. Also the probability mass function (pmf) is defined by

$$f(i) = \begin{cases} \dfrac{1}{n}, & i = 1,2,3,...,n, \\ 0 & elsewhere. \end{cases}$$

It implies that for the $\chi'$-chromatic mean and the $\chi'$-chromatic variance of $K_n$ we have

$$\mu_{\chi'}(K_n) = \dfrac{n+1}{2} \text{ and } \sigma_{\chi'}^2(K_n) = \dfrac{n^2-1}{12}.$$

By symmetry argument it also follows that:

$$\mu_{\chi^+}(K_n) = \dfrac{n+1}{2} \text{ and } \sigma_{\chi^+}^2(K_n) = \dfrac{n^2-1}{12}. \qquad \square$$

## 2.6.2 Braided Complete Graphs

Figure 3.3 depicts $G = K_7 \oplus_3 K_5$ where $K_n \oplus_l K_m$ means that $l \leq m$ and $l \leq n$ and $K_n \cap K_m = K_l$. We say that the two complete graphs are 1-braided. It is easy to see that the operation is commutative for two 1-braided complete graphs. For three or more 1-braided complete graphs we consider the graphs to be an ordered string of 1-braided complete graphs with the further condition that, if for $K_n \oplus_{l_1} K_m \oplus_{l_2} K_t$ we have that $K_n \cap K_m = K_{l_1}$ and $K_m \cap K_t = K_{l_2}$ then, $K_{l_1} \cap K_{l_2} = \emptyset$.

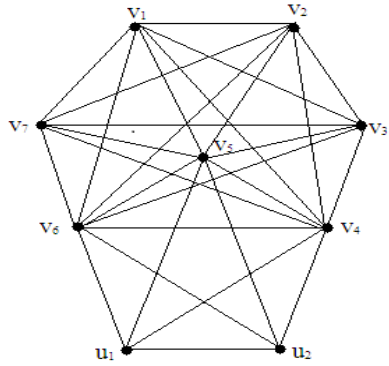

Fig 3.3

The probability mass function for the minimum chromatic sum $\chi'$ is

$$f(i) = \begin{cases} \dfrac{2}{n+m-l} = \dfrac{2}{9}, & i = 1,2, \\ \dfrac{1}{n+m-l} = \dfrac{1}{9}, & i = 3,4,,5,6,7, \\ 0, & \text{elsewhere}. \end{cases}$$

$$\therefore \mu_{\chi'}(G) = 1.\dfrac{2}{n+m-l} + 2.\dfrac{2}{n+m-l} + \sum_{i=3}^{7} i.\dfrac{1}{n+m-l} = \dfrac{6}{9} + \dfrac{25}{9} = \dfrac{31}{9}.$$

The probability mass function for the minimum chromatic sum $\chi^+$ is

$$f(i) = \begin{cases} \dfrac{2}{n+m-l} = \dfrac{2}{9}, & i = 6,7, \\ \dfrac{1}{n+m-l} = \dfrac{1}{9}, & i = 1,2,3,4,5, \\ 0, & \text{elsewhere}. \end{cases}$$

$$\therefore \mu_{\chi^+}(G) = 6 \cdot \frac{2}{n+m-l} + 7 \cdot \frac{2}{n+m-l} + \sum_{i=1}^{5} i \cdot \frac{1}{n+m-l} = \frac{26}{9} + \frac{15}{9} = \frac{41}{9}.$$

Furthermore, $\sigma^2_{\chi}(G) = 1^2 \cdot \frac{2}{n+m-l} + 2^2 \cdot \frac{2}{n+m-l} + \sum_{i=3}^{7} i^2 \cdot \frac{1}{n+m-l} - \left(\frac{31}{9}\right)^2 = \frac{145}{9} - \frac{961}{81} = \frac{344}{81}.$

and

$$\sigma^2_{\chi^+}(G) = 6^2 \cdot \frac{2}{n+m-l} + 7^2 \cdot \frac{2}{n+m-l} + \sum_{i=1}^{5} i^2 \cdot \frac{1}{n+m-l} - \left(\frac{41}{9}\right)^2 = \frac{255}{9} - \frac{1681}{81} = \frac{614}{81}.$$

**Lemma 2.6.3**

For two 1-braided complete graphs we have $\mu_{\chi'}(G) = \dfrac{2(n+1)(m-l) + (n-m+l)(n-m+l+1)}{2(n+m-l)}$

and $\mu_{\chi^+}(G) = \dfrac{(n-l)(n-l+1) + 4l(n-l) + 2l(l+1)}{2(n+m-l)}$

**Proof:**

**Part 1:** Without loss of generality assume $m \leq n$. The probability mass function for minimum chromatic sum $\chi'$ is

$$f(i) = \begin{cases} \dfrac{2}{n+m-l}, & i = 1,2,\ldots,m-l, \\ \dfrac{1}{n+m-l}, & i = (m-l)+1, (m-l)+2, \ldots, n, \\ 0, & \text{elsewhere.} \end{cases}$$

Therefore, $\mu_{\chi'}(G) = 1 \cdot \dfrac{2}{n+m-l} + 2 \dfrac{2}{n+m-l} + \ldots + (m-l)\dfrac{2}{n+m-l} +$

$\left\{ [(m-l)+1]\dfrac{1}{n+m-l} + [(m-l)+2]\dfrac{1}{n+m-l} + \ldots + n\dfrac{1}{n+m-l} \right\}$

$= \dfrac{2}{n+m-l}[1+2+3+\ldots+(m-l)] + \dfrac{1}{n+m-l}[(m-l)+1+(m-l)+2+\ldots+n]$

$= \dfrac{2}{n+m-l}\left[\dfrac{(m-l)(m-l+1)}{2}\right] + \dfrac{1}{n+m-l}[(n-m+l)(m-l)+1+2+3+\ldots+(n-m+l)]$

$= \dfrac{(m-l)(m-l+1)}{n+m-l} + \dfrac{1}{n+m-l}[(n-m+1)(m-l) + \dfrac{(n-m+l)(n-m+l+1)}{2}]$

$$= \frac{1}{n+m-l}\left\{(m-l)(m-l+1)+(n-m+l)(m-l)+\frac{(n-m+l)(n-m+l+1)}{2}\right\}.$$

$$\therefore \mu_{\chi^-}(G) = \frac{2(m-l)(n+1)+(n-m+l)(n-m+l+1)}{2(n+m-l)}.$$

**Part 2:** Without loss of generality assume $m \leq n$. The probability mass function (pmf) for maximum chromatic sum $\chi^+$ is

$$f(i) = \begin{cases} \dfrac{2}{n+m-l}, & i = (n-l+1), (n-l+2), (n-l+3), \ldots, ((n-l)+l), \\ \dfrac{1}{n+m-l}, & i = 1, 2, 3, \ldots, n-l, \\ 0, & \text{elsewhere}. \end{cases}$$

Therefore
$$\mu_{\chi^+}(G) = \left\{1 \cdot \frac{1}{n+m-l} + 2 \cdot \frac{1}{n+m-l} + \ldots + (n-l) \cdot \frac{1}{n+m-l}\right\} +$$
$$\left\{(n-l+1) \cdot \frac{2}{n+m-l} + (n-l+2) \cdot \frac{2}{n+m-l} + \ldots + ((n-l)+l) \cdot \frac{2}{n+m-l}\right\}$$

$$= \frac{1}{n+m-l}\{1+2+3+\ldots+(n-l)\} + \frac{2}{n+m-l}\{(n-l+1)+(n-l+2)+\ldots+((n-l)+l)\}$$

$$= \frac{1}{n+m-l}\left\{\frac{(n-l)(n-l+1)}{2}\right\} + \frac{2}{n+m-l}\{l(n-l)+1+2+3+\ldots+l\}$$

$$= \frac{1}{n+m-l}\left\{\frac{(n-l)(n-l+1)}{2}\right\} + \frac{2}{n+m-l}\left\{l(n-l)+\frac{l(l+1)}{2}\right\}$$

$$= \frac{1}{n+m-l}\left\{\frac{(n-l)(n-l+1)}{2}\right\} + \frac{2}{n+m-l}\left\{\frac{2l(n-l)+l(l+1)}{2}\right\}$$

$$\therefore \mu_{\chi^+}(G) = \frac{(n-l)(n-l+1)+4l(n-l)+2l(l+1)}{2(n+m-l)}. \qquad \square$$

### 2.6.5 Application to $q$-Jaco Graph for: $f(x) = x^2$

Informally stated, it is noted that $q$-Jaco graphs have the structural property that they are typically a string of intersecting complete graphs with the complete graphs following in non-decreasing order as $n$ increases. In this section the introductory results are applied to the $q$-Jaco graph for $f(x) = x^2$.

The corresponding colouring sum weights for $\chi'(J_i(f(x)))$, $\chi^+(J_i(f(x)))$, are depicted in the table 2.

Table 2

| $\phi(v_i) \to i \in N$ | Minimum Chromatic Sum Colour Weights | Maximum Chromatic Sum Colour Weights |
|---|---|---|
| 1 | $\theta(c_1) = 1$ | $\theta(c_1) = 1$ |
| *2 | $\theta(c_1) = 1$, $\theta(c_2) = 1$ | $\theta(c_1) = 1$, $\theta(c_2) = 1$ |
| 3 | $\theta(c_1) = 2$, $\theta(c_2) = 1$ | $\theta(c_1) = 1$, $\theta(c_2) = 2$ |
| 4 | $\theta(c_1) = 2$, $\theta(c_2) = 1$, $\theta(c_3) = 1$ | $\theta(c_1) = 1$, $\theta(c_2) = 1$, $\theta(c_3) = 2$ |
| *5 | $\theta(c_1) = 2$, $\theta(c_2) = 1$, $\theta(c_3) = 1$, $\theta(c_4) = 1$ | $\theta(c_1) = 1$, $\theta(c_2) = 1$, $\theta(c_3) = 1$, $\theta(c_4) = 2$ |
| 6 | $\theta(c_1) = 2$, $\theta(c_2) = 2$, $\theta(c_3) = 1$, $\theta(c_4) = 1$ | $\theta(c_1) = 1$, $\theta(c_2) = 1$, $\theta(c_3) = 2$, $\theta(c_4) = 2$ |
| 7 | $\theta(c_1) = 2$, $\theta(c_2) = 2$, $\theta(c_3) = 1$, $\theta(c_4) = 1$, $\theta(c_5) = 1$ | $\theta(c_1) = 1$, $\theta(c_2) = 1$, $\theta(c_3) = 1$, $\theta(c_4) = 2$, $\theta(c_5) = 2$ |
| 8 | $\theta(c_1) = 2$, $\theta(c_2) = 2$, $\theta(c_3) = 1$, $\theta(c_4) = 1$, $\theta(c_5) = 1$, $\theta(c_6) = 1$ | $\theta(c_1) = 1$, $\theta(c_2) = 1$, $\theta(c_3) = 1$, $\theta(c_4) = 1$, $\theta(c_5) = 2$, $\theta(c_6) = 2$ |
| 9 | $\theta(c_1) = 2$, $\theta(c_2) = 2$, $\theta(c_3) = 1, \ldots, \theta(c_7) = 1$ | $\theta(c_1) = 1, \ldots, \theta(c_5) = 1$, $\theta(c_6) = 2$, $\theta(c_7) = 2$ |
| 10 | $\theta(c_1) = 2$, $\theta(c_2) = 2$, $\theta(c_3) = 1, \ldots, \theta(c_8) = 1$ | $\theta(c_1) = 1, \ldots, \theta(c_6) = 1$, $\theta(c_7) = 2$, $\theta(c_8) = 2$ |
| *11 | $\theta(c_1) = 2$, $\theta(c_2) = 2$, $\theta(c_3) = 1, \ldots, \theta(c_9) = 1$ | $\theta(c_1) = 1, \ldots, \theta(c_7) = 1$, $\theta(c_8) = 2$, $\theta(c_9) = 2$ |
| 12 | $\theta(c_1) = 3$, $\theta(c_2) = 2$, $\theta(c_3) = 1, \ldots, \theta(c_9) = 1$ | $\theta(c_1) = 1, \ldots, \theta(c_7) = 1$, $\theta(c_8) = 2$, $\theta(c_9) = 3$ |
| 13 | $\theta(c_1) = 3$, $\theta(c_2) = 2$, | $\theta(c_1) = 1, \ldots, \theta(c_8) = 1$, |

|    | $\theta(c_3)=1,\ldots,\theta(c_{10})=1$ | $\theta(c_9)=2,\theta(c_{10})=3$ |
|----|---|---|
| 14 | $\theta(c_1)=3,\theta(c_2)=2,$ $\theta(c_3)=1,\ldots,\theta(c_{11})=1$ | $\theta(c_1)=1,\ldots,\theta(c_9)=1,$ $\theta(c_{10})=2,\theta(c_{11})=3$ |
| 15 | $\theta(c_1)=3,\theta(c_2)=2,$ $\theta(c_3)=1,\ldots,\theta(c_{12})=1$ | $\theta(c_1)=1,\ldots,\theta(c_{10})=1,$ $\theta(c_{11})=2,\theta(c_{12})=3$ |
| 16 | $\theta(c_1)=3,\theta(c_2)=2,$ $\theta(c_3)=1,\ldots,\theta(c_{13})=1$ | $\theta(c_1)=1,\ldots,\theta(c_{11})=1,$ $\theta(c_{12})=2,\theta(c_{13})=3$ |
| 17 | $\theta(c_1)=3,\theta(c_2)=2,$ $\theta(c_3)=1,\ldots,\theta(c_{14})=1$ | $\theta(c_1)=1,\ldots,\theta(c_{12})=1,$ $\theta(c_{13})=2,\theta(c_{14})=3$ |
| *18 | $\theta(c_1)=3,\theta(c_2)=2,$ $\theta(c_3)=1,\ldots,\theta(c_{15})=1$ | $\theta(c_1)=1,\ldots,\theta(c_{13})=1,$ $\theta(c_{14})=2,\theta(c_{15})=3$ |
| 19 | $\theta(c_1)=3,\theta(c_2)=2,\theta(c_3)=2,$ $\theta(c_4)=1,\ldots,\theta(c_{15})=1$ | $\theta(c_1)=1,\ldots,\theta(c_{12})=1,$ $\theta(c_{13})=2,\theta(c_{14})=2,\theta(c_{15})=3$ |
| 20 | $\theta(c_1)=3,\theta(c_2)=2,\theta(c_3)=2,$ $\theta(c_4)=1,\ldots,\theta(c_{16})=1$ | $\theta(c_1)=1,\ldots,\theta(c_{13})=1,$ $\theta(c_{14})=2,\theta(c_{15})=2,\theta(c_{16})=3$ |

Table 3 below depicts the values $\chi'(J_i(f(x)))$, $\chi^+(J_i(f(x)))$, $\mu_{\chi'}(J_i(f(x)))$, $\mu_{\chi^+}(J_i(f(x)))$, $\sigma^2_{\chi'}(J_i(f(x)))$, $\sigma^2_{\chi^+}(J_i(f(x)))$ for $1 \leq i \leq 20$.

Table 3

| $\phi(v_i) \to i \in N$ | $\chi'(J_i(f))$ | $\chi^+(J_i(f))$ | $\mu_{\chi'}(J_i(f))$ | $\mu_{\chi^+}(J_i(f))$ | $\sigma^2_{\chi'}(J_i(f))$ | $\sigma^2_{\chi^+}(J_i(f))$ |
|---|---|---|---|---|---|---|
| 1 | 1 | 1 | 1 | 1 | 0 | 0 |
| *2 | 3 | 3 | $3/2$ | $3/2$ | $1/4$ | $1/4$ |
| 3 | 4 | 5 | $4/3$ | $5/3$ | $2/9$ | $2/9$ |
| 4 | 7 | 9 | $7/4$ | $9/4$ | $11/16$ | $11/16$ |
| *5 | 11 | 14 | $11/5$ | $14/5$ | $34/25$ | $34/25$ |
| 6 | 13 | 17 | $13/6$ | $17/6$ | $41/36$ | $41/36$ |
| 7 | 18 | 24 | $18/7$ | $24/7$ | $96/49$ | $96/49$ |
| 8 | 24 | 32 | $24/8$ | $32/8$ | $192/64$ | $192/64$ |
| 9 | 31 | 41 | $31/9$ | $41/9$ | $344/81$ | $344/81$ |
| 10 | 39 | 51 | $39/10$ | $51/10$ | $469/100$ | $469/100$ |
| *11 | 48 | 62 | $48/11$ | $62/11$ | $886/121$ | $886/121$ |
| 12 | 49 | 71 | $49/12$ | $71/12$ | $1091/144$ | $1091/144$ |
| 13 | 59 | 84 | $59/13$ | $84/13$ | $1602/169$ | $1602/169$ |
| 14 | 70 | 98 | $70/14$ | $98/14$ | $2268/196$ | $2268/196$ |
| 15 | 82 | 113 | $82/15$ | $113/15$ | $3116/225$ | $3116/225$ |
| 16 | 95 | 129 | $95/16$ | $129/16$ | $4175/256$ | $4175/256$ |
| 17 | 104 | 146 | $109/17$ | $146/17$ | $5476/289$ | $5476/289$ |
| *18 | 119 | 164 | $124/18$ | $164/18$ | $7852/324$ | $7852/324$ |
| 19 | 122 | 177 | $127/19$ | $177/19$ | $7716/361$ | $7716/361$ |
| 20 | 138 | 197 | $143/20$ | $197/20$ | $9771/20$ | $9771/20$ |

We observe that the colour weights undergo well-defined change as the order of the $q$-Jaco graph increases from $i$ to $i + 1$. We present the next result in this respect.

**Conjecture 2.6.4**

For the $q$-Jaco graph $J_n(f(x)), n \in N$ and for $f(x) = x^2$, we have $\sigma^2_{\chi}(J_i(f(x))) = \sigma^2_{\chi^+}(J_i(f(x)))$.

**Theorem 2.6.5**

When the order of the $q$-Jaco graph with $f(x) = x^2$ increases from $i$ to $i + 1$, the minimum sum colouring changes by:

  (i)   Either a new colour with minimum subscript say, $c_k$ is added with $\theta(c_k) = 1$ or,

  (ii)  Exactly one of the colours in $J_i(f(x))$ increases by count 1.

**Proof:**

If the degree of the prime Jaconian vertex $v_p$ of $J_i(f(x))$ has not reached its defined maximum, $p^2$ it implies that the corresponding Hope graph increases in order by +1 as $J_i(f(x))$ expands to $J_{i+1}(f(x))$. Clearly, in terms of the definition of a proper colouring a new colour with minimum subscript say, $c_k$ is added with $\theta(c_k) = 1$.

If the degree of the prime Jaconian vertex $v_p$ of $J_i(f(x))$ has reached its defined maximum, $p^2$ it implies that as $J_i(f(x))$ expands to $J_{i+1}(f(x))$ it may be coloured with an existing colour with minimum subscript allowed by a proper colouring of the Hope graph. □

## 3. CONCLUSION

Many characteristics of linear Jaco graph were found in this paper. The concept of quadratic Jaco graph has been introduced. We also derived some results and characteristics of quadratic Jaco graph. Similar results could be extended to polynomial Jaco graph too. The aforesaid will be the subject of further research.

**REFERENCES**

[1] J.A. Bondy and U. S. R. Murty, **Graph Theory with Applications**, Macmillan Press, London (1976).

[2] G. Chartrand and P. Zang, **Chromatic graph theory**, CRC Press, 2009.

[3] T.R. Jensen and B. Toft, **Graph coloring problems**, John Wiley & Sons, 1995.

[4] J. Kok, P. Fisher, B. Wilkens, M. Mabula and V. Mukungunungwa, *Characteristics of finite Jaco graphs, $J_n(1), n \in N$,* arXiv: 1404.0484v1 [math.CO], 2 April 2014.

[5] J. Kok, P. Fisher, B. Wilkens, M. Mabula and V. Mukungunungwa, "*Characteristics of Jaco Graphs, $J_\infty(a), n \in N$,* arXiv: 1404.0484v1 [math.CO], 7 April 2014.

[6] J. Kok, C. Susanth and Sunny Joseph Kalayathankal, *A Study on Linear Jaco Graphs, Journal of Informatics and Mathematical Sciences*", Vol.7, No. 2, pp. 69-80, 2015.

[7] J. Kok, *Linear Jaco Graphs: A Critical Review*, Journal of Informatics and Mathematical Sciences, Vol. 8, No. 2, pp. 67-103, 2016.

[8] J. Kok, N.K. Sudev and K.P. Chithra, *General colouring sums of graphs*, Cogent Mathematics, (2016), **3**, 1140002.

[9] N.K. Sudev, S. Satheesh, K.P. Chithra and J. Kok, *On Certain Colouring Parameters of Graphs*, Communicated.